\documentclass[twocolumn]{autart}    

\usepackage{graphicx}          
\graphicspath{{fig/}}
\newcommand\diff{\,\mathrm{d}}
\usepackage{color}
\usepackage[cmyk]{xcolor}
\usepackage{amsfonts}
\usepackage{graphicx,psfrag}
\usepackage{amsmath, bm,amssymb}
\usepackage{subfig}

\allowdisplaybreaks[4]

\newtheorem{theorem}{Theorem}

\newtheorem{remark}{Remark}
\begin{document}

\begin{frontmatter}

\title{ Robust stabilization of  $2 \times 2$  first-order hyperbolic PDEs  with uncertain  input delay\thanksref{t1} \thanksref{t2}} 

\thanks[t1]{This paper was not presented at at any conference.}
\thanks[t2]{This work was supported by the National Natural Science Foundation of China (62173084, 61773112), the Project of Science and Technology Commission of Shanghai Municipality, China (23ZR1401800, 22JC1401403).}

\author[DHU]{Jing Zhang,} \ead{zhangjing@mail.dhu.edu.cn} 
\author[DHU,Textile]{Jie Qi\thanksref{cor}} \ead{jieqi@dhu.edu.cn}           
\address[DHU]{College of Information Science and Technology, Donghua University, Shanghai, 201620, China} 
\address[Textile]{Engineering Research Center of Digitized Textile and Fashion Technology Ministry of Education, Donghua University, Shanghai, 201620, China}  
\thanks[cor]{Corresponding author.}
          
\begin{abstract}                         
 A backstepping-based compensator design is developed for a system of $2\times2$ first-order linear hyperbolic partial differential equations (PDE) in the presence of an uncertain long input delay at boundary. We introduce a transport PDE to represent the delayed input, which leads to three coupled first-order hyperbolic PDEs. A novel backstepping transformation, composed of two Volterra transformations and an affine Volterra transformation, is introduced for the predictive control design. The resulting kernel equations from the affine Volterra transformation are two coupled first-order PDEs and each with two boundary conditions, which brings challenges to the well-posedness analysis. We solve the challenge by using the method of characteristics and the successive approximation.  To analyze the sensitivity of the closed-loop system to uncertain input delay, we introduce a neutral system which captures the control effect resulted from the delay uncertainty. It is proved that the proposed control is robust to small delay variations.  Numerical examples illustrate the performance of the proposed compensator. 
\end{abstract}
\end{frontmatter}

\section{Introduction}
 The compensator design for distributed parameter systems (DPS) with delay, which is molded by PDEs, has been given much concern and study  \cite{xu2006stabilization,lhachemi2021robustness,krstic2008backstepping}. 
 Delays widely exist in industrial process fields and common phenomena, yet even minor delays may destabilize the system or cause accidents. Hence, the compensation control of delay is an active topic. The main methods of spectrum decomposition \cite{xu2006stabilization,braik2022well}, Lyapunov function technique \cite{prieur2012iss}, fuzzy control \cite{ding2021fuzzy}, Hilbert space theory \cite{khusainov2015strong}, and backstepping approach \cite{krstic2009delay} are utilized to stabilize such systems with delay.  

 As well known, unbounded control operators based on the backstepping method have been considered for both  hyperbolic equations \cite{zhang2021compensation},  heat equations \cite{krstic2009control,qi2020compensation}, wave equations, and some high-dimensional PDE systems\cite{qi2019control}, where delays are compensated. Not only delay compensation design, but also a series of problems in DPS  can be solved by backstepping method: motion
 planning  \cite{freudenthaler2020pde},  adaptive controller design \cite{anfinsen2017adaptive,wang2021delay} and output regulation problem \cite{deutscher2015backstepping,zhang2022output}, etc.
 Among them is a class of PDE systems, $2\times 2$ first-order hyperbolic PDEs, which can be used to describe various industrial processes and phenomena, such as traffic flow \cite{colombo20022}, 
 hydraulics and river dynamics  \cite{dos2008boundary,bastin2011boundary}, pipeline networks  \cite{aamo2015leak} and so on. \cite{vazquez2011backstepping} stabilizes   $2\times 2$ linear system by backstepping method, which is the foundation of subsequent works,  and further, extending to quasi-linear system in \cite{coron2013lloca}.  The more general stabilization problem of a class of   $m+1$  and $m+n$  hyperbolic PDEs system  are addressed in \cite{di2013stabilization,hu2016control}. The disturbance rejection and observer design for hyperbolic systems have been focused on in \cite{tang2014active,anfinsen2014disturbance} and \cite{irscheid2021observer}, respectively.  An event-triggered boundary control is proposed in \cite{wang2022delay}, where the sensor delay is compensated.
  
 The delay-compensated control design is more challenging for uncertain coupled hyperbolic systems.  One way for robustly stabilizing a $2\times 2$ linear hyperbolic system against possible actuator delay by preserving a small amount of the proximal reflection is  presented in \cite{auriol2018delay}. Similarly, a $m+n$ hyperbolic PDEs system is considered in \cite{auriol2019explicit}. 
An underactuated cascade network of $n$ systems of two heterodirectional hyperbolic PDEs is considered in \cite{auriol2020outputtwo}, which utilizes a two-step integral transformation to move the in-domain coupling terms to the boundary.    
In addition, the method of combining with a low-pass filter and a controller based on prediction is proposed  in \cite{auriol2022robust}, for solving the robust stabilization problem of an underactuated network of two subsystems of heterodirectional $n+m$ hyperbolic system. The method for recursively handling multiple cascaded subsystems is also applied to stabilize systems that concatenate  multiple hyperbolic PDEs and an ODE  \cite{auriol2021stabilization,redaud2021output}.  

 Another approach to address uncertain delays is to design adaptive control laws using existing  stability results. For instance, the $2\times2$ hyperbolic system with input delay and sensor delay is equivalent to a delay-free system in \cite{anfinsen2018adaptive} through a series of transformations; then the existing adaptive control law can be employed.

 This paper provides a distinctive way of stabilization of the $2 \times 2$ hyperbolic system subject to an uncertain input delay. First, the exact system is taken into account. A transport PDE is introduced to represent the delayed input, which results in three coupled cascade first-order hyperbolic PDEs.  Unlike the existing references which apply multi-step transformation in the control design, such as \cite{auriol2018delay}, we propose a more direct way to solve the delay-compensation problem by introducing a single-step backstepping transformation composed of two Volterra transformations  and an affine Volterra transformation. The affine Volterra transformation includes  both the integration of the states and that of the delayed actuator state governed by the transport PDE. The resulting kernel equations of the affine Volterra transformation are coupled both in the domain equations and in the boundary conditions. Besides, there are two boundary conditions for each kernel, which leads to each kernel being defined in multiple subregions. Due to these two facts, it is challenging to prove the well-posedness of the kernel equations. The challenge is solved by combining the method of characteristics and the successive approximation.    
  
 In line with the transformation, the control law consists of the feedback of the states and the feedback of historical control input. The finite-time stability of closed-loop systems with delay-compensated control is proved.  Furthermore, we analyze the robustness of the delay-compensated control for the uncertain input delay through the equivalent neutral system.  Finally, compared with the controller without delay compensation, the control effect of the proposed compensator is illustrated. 

The contribution of this paper is summarized as: 
\begin{enumerate}
	\item  Compared to the results in \cite{auriol2018delay,auriol2021stabilization,redaud2021output} about robust stabilization of hyperbolic system with small actuator delay, we propose an arbitrary long delay-compensated control design method and analyze the robustness of the controller with small delay variation. 
	
	\item  We develop a novel single-step backstepping transformation that is more straightforward than the multi-step transformation used in  \cite{auriol2020outputtwo,anfinsen2018adaptive} because it avoids multi-step inverse transformation to reach the final form of the control law. The resulting kernel equations appear more complicated form than that of  \cite{vazquez2011backstepping} due to each kernel equation having two boundary conditions. We combine the method of characteristics and the successive approximation to prove that the kernel equations, where each kernel function is defined in multiple subregions, are well-posed. 
\end{enumerate}

 This paper is organized as follows.   Section \ref{sec:PDE-system} gives the problem formulation.  The backstepping-based  compensation controller is established in Section \ref{sec:control}. In Section \ref{sec:kernel}, well-posedness of the kernels is discussed, which is one of the main technical contribution of this paper. The finite-time stability of the closed-loop system that is subjected to exact input delay is summarized in section \ref{sec:stability}.  In Section \ref{sec:robust}, we prove that the proposed control law still stabilizes the hyperbolic system if there is small variation in the input delay, which also is one of the technical contributions of this paper. The simulation results are presented in Section \ref{sec:simulation} to validate the effectiveness of the proposed compensator. 
 The paper ends with concluding remarks and a discussion of the future work in Section \ref{sec:conclusion}. 
\section{Problem formulation}
\label{sec:PDE-system}
Consider the following  $2 \times 2$ first-order linear hyperbolic PDEs with input delay $\tau\in R^+$,  
\begin{subequations}\label{eq:system}
	\begin{align}
		\partial_t  u_1 (x ,t)  =&~  -\varepsilon_{1}(x ) 	\partial_x u_1(x ,t) +c_1(x ) u_2(x ,t) , \label{eq:main-u1}\\  
		\partial_t  u_2 (x ,t)  =&~    \varepsilon_{2}(x ) \partial_x u_2(x ,t) +c_2(x ) u_1(x ,t) , \label{eq:main-u2}\\
		u_1(0,t)  = &~    q u_2(0,t),\label{eq:bnd-u1}\\
		u_2(1,t)  =&~      U(t-\tau),\label{eq:bnd-u2}\\
		U(t-h)  =&~ 0, \quad h\in[0,\tau], 
	\end{align}
\end{subequations}
for $(x,t)  \in [0,1] \times \mathbb{R} ^+$,   with the initial conditions,
\begin{subequations}
	\begin{align}
		u_1(x,0)=u_{10}(x)\in\mathbb{C}([0,1]),\\
		u_2(x,0)=u_{20}(x)\in\mathbb{C}([0,1]),
	\end{align}
\end{subequations}
and     where the coefficients $q\in \mathbb{R}^+$, $\varepsilon_{1}(x)$, $\varepsilon_{2}(x), c_1(x)$,  $c_2(x) \in \mathbb{C}([0,1])$ and transport velocities $\varepsilon_{1}(x),\,\varepsilon_2(x)\textgreater 0$. Here,  $\partial_x  u(x ,t)= {\partial u}/{\partial x }$ and  $\partial_t u(x ,t)= {\partial    u }/{\partial t}$ are represented in a compact form of partial derivative.

The control objective of this paper is to design a compensated controller to ensure that the closed-loop system is exponentially stable, in particular, reaching the zero-equilibrium state in finite time.   
Based on the relation between the first-order hyperbolic PDE and delay equations \cite{karafyllis2014relation}, a transport PDE is introduced to express the input delay, so the system \eqref{eq:system} is rewritten as follows: 
\begin{subequations}\label{eq:coupled_system}
	\begin{align}
		\partial_t  u_1 (x ,t) &=  -\varepsilon_{1}(x )\partial_x  u_1(x ,t) +c_1(x ) u_2(x ,t) ,  \label{eq:main-coupled-u1}\\  
		\partial_t  u_2 (x ,t) &=   \varepsilon_{2}(x )\partial_x  u_2(x ,t) +c_2(x ) u_1(x ,t) ,  \label{eq:main-coupled-u2}\\
		u_1(0,t) &=    q u_2(0,t),\label{eq:bnd-coupled-u1}\\
		u_2(1,t) &=   v(0,t),\label{eq:bnd-coupled-u2}\\
		\tau \partial_{t} v(x ,t)&=\partial_{x } v(x ,t), \\
		v(1,t)&=U(t), \label{eq:coupled_system_v_bud}
	\end{align}
\end{subequations}
with the initial condition $v(x,0)=v_0(x)\in\mathbb{C}([0,1])$.
\section{Controller design} \label{sec:control}
\subsection{Backstepping control design} \label{sec_sub:back}
Following the backstepping approach \cite{krstic2008boundary_book}, we introduce  two Volterra transformations and an affine Volterra transformation integral transformation as follows:
\begin{subequations}\label{eq:map}
	\begin{align}
		w_1(x ,t)=&~   u_1(x ,t) - \int_{0}^{x }  K^{11}(x ,y) u_1(y,t) \diff y\nonumber\\
		&~-\int_{0}^{x } K^{12}(x ,y) u_2(y,t) \diff y, \label{eq:map1}\\
		w_2(x ,t)=&~  u_2(x ,t) - \int_{0}^{x } \!K^{21}(x ,y) u_1(y,t) \diff y\nonumber\\
		&~-\int_{0}^{x } K^{22}(x ,y) u_2(y,t) \diff y, \label{eq:map2}
		\\
		z(x ,t) 
		=&~v(x ,t)-\int_{0}^{x } p (x -y) v(y,t) \diff y\label{eq:map3}\\
		-  \int_{0}^{1}&  \alpha_1 (x ,y) u_1(y,t) \diff y - \int_{0}^{1}  \alpha_2(x ,y) u_2(y,t) \diff y,\nonumber 
	\end{align}
\end{subequations}
where the  kernel  functions $K^{ij}(x,y)$, $i,j=1,2$, are  to be defined on the triangular domain $\mathcal T_1 =\left\{(x ,y)\in \mathbb{R}^2:0\leq y  \leq x \leq 1\right\}$,  the  kernel  functions $\alpha_1(x,y)$, $\alpha_2(x,y)$ are  to be defined on the rectangle domain $\mathcal T_2 =\left\{(x ,y)\in \mathbb{R}^2:0\leq x,\,y \leq 1\right\}$, and the kernel function $p(x)\in L^2([0,1])$. 

Then, we find sufficient conditions on these functions, such that transformations \eqref{eq:map} map system \eqref{eq:coupled_system} to the following system
\begin{subequations}\label{eq:target}
	\begin{align}
		\partial_t  w_1 (x ,t) &=  -\varepsilon_{1}(x )\partial_x w_1(x ,t)  , \label{eq:main-coupled-w1}\\  
		\partial_t  w_2 (x ,t) &= \varepsilon_{2}(x )\partial_x w_2(x ,t)   , \label{eq:main-coupled-w2}\\
		w_1(0,t) &= q w_2(0,t),\label{eq:bnd-coupled-w1}\\
		w_2(1,t) &=  z(0,t),\label{eq:bnd-coupled-w2}\\
		\tau \partial_{t} z(x ,t)&=\partial_{x } z(x ,t), \label{eq:target_v} \\
		z(1,t)&=0.\label{eq:target_v_bud}
	\end{align}
\end{subequations}
Combined with system \eqref{eq:coupled_system}, differentiate the mapping \eqref{eq:map} with respect to space and time respectively, and then substitute it into the target system \eqref{eq:target},  one obtains that the kernels $K^{ij} $, $i,j=1,2$, satisfy a set of hyperbolic PDEs given in \cite{coron2013lloca},
and the kernels $\alpha_1$, $\alpha_2$ satisfy:
\begin{subequations}\label{eq:alpha_main}
	\begin{align}
		\frac{1}{\tau}\partial_x  \alpha_1  -  \varepsilon_{1}(y)& \partial_y \alpha_1  =\varepsilon_{1}'(y) \alpha_1(x ,y)   +c_2(y)\alpha_2 ,  \label{eq:alpha1}\\
		\frac{1}{\tau}\partial_x  \alpha_2 +  \varepsilon_{2}(y)& \partial_y \alpha_2 =-\varepsilon_{2}'(y) \alpha_2  +c_1(y)\alpha_1 , \label{eq:alpha2}	 
	\end{align}
\end{subequations}
with the boundary conditions
\begin{subequations}\label{eq:alpha_bud}
	\begin{align}
		\alpha_1(x ,1)=&~0,\label{eq:alpha_bud2} \\
		\varepsilon_2(0)\alpha_2(x ,0)=&~ q\varepsilon_{1}(0) \alpha_1(x ,0),\label{eq:beta_bud1}\\
		\alpha_1(0,y)=&~K^{21}(1,y) ,\label{eq:alpha_bud1} \\
		\alpha_2(0,y)=&~K^{22}(1,y) . \label{eq:beta_bud2}
	\end{align}
\end{subequations}
After solving $\alpha_2(x, y)$, the kernel $p(x)$ has the following form
\begin{align}\label{eq:p}
	p(x)=\tau \varepsilon_{2}(1) \alpha_2(x
	,1).
\end{align}
After obtaining all the kernels, we discuss the form of compensator. Substituting boundary conditions \eqref{eq:coupled_system_v_bud} and \eqref{eq:target_v_bud} into transformation \eqref{eq:map3}, the controller is obtained in the following form:
\begin{align}\label{eq:U}
	U(t)= &~\int_{0}^{1} p (1-y) v(y,t) \diff y+\int_{0}^{1} \alpha_1 (1,y) u_1(y,t) \diff y\nonumber\\
	&~+\int_{0}^{1} \alpha_2(1,y) u_2(y,t) \diff y.
\end{align}
Then, using the method of
characteristics and variable substitution for \eqref{eq:U}, the controller with delay compensation can be written in the following form:
\begin{align}\label{eq:control}
	U(t)=& \int_{0}^{1} \alpha_1 (1,y) u_1(y,t) \diff y 
	+\int_{0}^{1} \alpha_2(1,y) u_2(y,t) \diff y\nonumber\\
	& 
	\int_{t-\tau}^{t}p\left(\frac{t-y}{\tau}\right) U(y) \diff y.
\end{align}

\section{Well-posedness of the kernel functions} \label{sec:kernel}
The well-posedness of the kernel functions $K^{ij}$, $i,j=1,2$ have given in \cite{coron2013lloca}. Here we discuss the well-posedness of $\alpha_1$ and $\alpha_2$, because they have two boundary conditions which make them different from kernels $K^{ij}$ or $L^{ij}$ with only one boundary condition in a triangle domain. Define a general form of  equations $\alpha_1$-$\alpha_2$   as follows:
\begin{subequations}\label{eq:ker_gen}
	\begin{align}
		\frac{\partial_x  F_1}{\tau}  - \varepsilon_{1}(y) \partial_y F  & = g_1(x,y) +\sum_{k=1}^{2}C_{1k}(x,y)F_k , \label{ker_gen_1} \\
		\frac{\partial_x  F_2}{\tau}  + \varepsilon_{2}(y) \partial_y F_2 & =\! g_2(x,y) +\sum_{k=1}^{2}C_{2k}(x,y)F_k , \label{ker_gen_2}
	\end{align}	 
\end{subequations}
evolving in the domain $\mathcal T_2$, with boundary conditions
\begin{subequations}\label{eq:ker_gen_bud}
	\begin{align}
		F_1(x ,1)&=0,\label{ker_gen_bud1}\\
		F_1(0,y)&=h_1(y),\label{ker_gen_ini1}\\
		F_2(x ,0)&= q_1(x) F_1(x ,0),\label{ker_gen_bud2}\\	
		F_2 (0,y)&=h_2(y).   \label{ker_gen_ini3}
	\end{align}	 
\end{subequations}
The well-posedness of equations \eqref{eq:ker_gen} and  \eqref{eq:ker_gen_bud}  is presented in the following theorem.
\begin{theorem} \label{theo:well} \normalfont 
	Consider the hyperbolic system \eqref{eq:ker_gen} with the boundary condition \eqref{eq:ker_gen_bud}. Under the assumptions  $q_1(x), h_i(x)\in \mathbb{C}([0,1])$, $\varepsilon_{i}(x)\in\mathbb{C}([0,1])$  with $\varepsilon_{i}(x)\textgreater 0$ and  $C_{ij}(x,y)\in \mathcal{T}_2$, $i,j=1,2$, there exists a unique  solution $F_k \in  \mathbb C(\mathcal T_2)$, $k =1, 2$.
\end{theorem}
The proof of this theorem is divided into two steps. First, we transform the coupled equations into integral equations using the method of characteristics, which is presented in Section \ref{sec:Integ}, and then solve it by successive approximation method, which is presented in Section \ref{sec:convergence}. 
\subsection{Integral equation} \label{sec:Integ}
First, define 
\begin{align}
	\phi_1(x)=\int_{0}^{x}\frac{1}{\varepsilon_{1}(\zeta)}\diff \zeta,\quad \phi_2(x)=\int_{0}^{x}\frac{1}{\varepsilon_{2}(\zeta)}\diff \zeta.
\end{align}
Note that since the propagation velocities  $\varepsilon_{1}(x)$ and $\varepsilon_{2}(x)$ are both positive, all $\phi_i(x)$, $i=1,2$ functions are monotonically increasing and invertible.
Also, it holds that $\phi_i, \phi_i^{-1} \in \mathbb{C}([0,1])$.
\begin{figure}[ht]    
	\centering
	\includegraphics[width=0.48\textwidth]{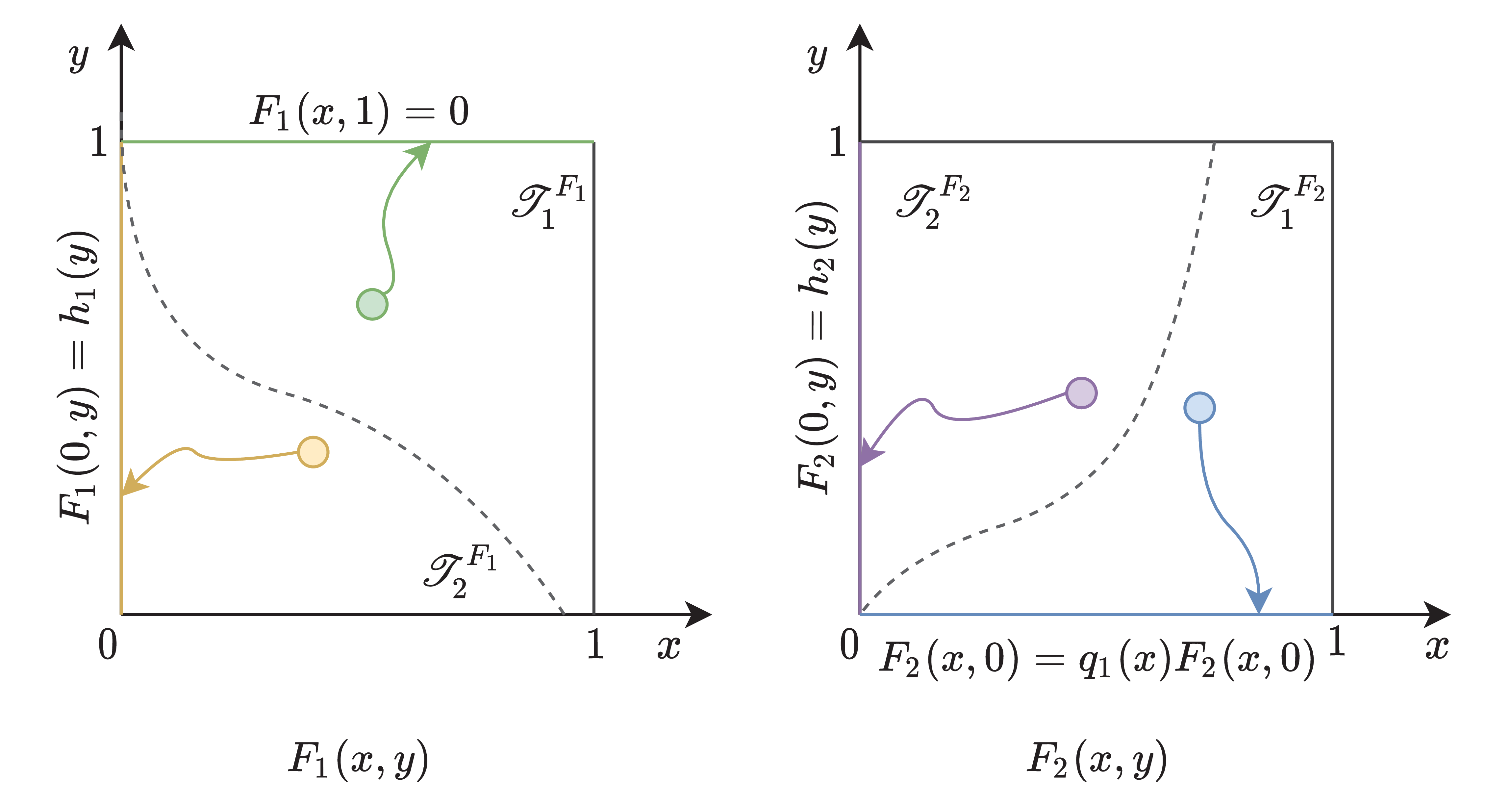}
	\caption{Diagram of characteristic curves. Different color curves represent characteristic curves affected by corresponding boundary conditions.}
	\label{fig:line}
\end{figure}
The characteristic curves of $F_1(x,y)$ and $F_2(x,y)$ are depicted in Fig. \ref{fig:line}. According to the distribution of the characteristic curves,  $F_1(x,y)$ is defined in two regions $\mathcal{T}^{F_1}_1$ and  $\mathcal{T}^{F_1}_2$,  and  $F_2(x,y)$ in $\mathcal{T}^{F_2}_1$ and $\mathcal{T}^{F_2}_2$, which are shown in Fig. \ref{fig:line}. It is means that, for instance, when solving the characteristic \eqref{eq:F1_c_x1}--\eqref{eq:F1_c_y1} backward from a given point in $\mathcal{T}^{F_1}_1$, one   ``hits'' the boundary depicted in green on Fig. \ref{fig:line}. 

(1) Characteristics of the $F_1(x,y)$:
For $(x,y)\in \mathcal{T}_2$, define the characteristic curves as follows, along which the system \eqref{ker_gen_1} becomes a family of ordinary differential equations, 
\begin{subequations}
	\begin{align}
		x_{11}(x,y,s)=&~x +\frac{s-\phi_1(1)+\phi_1(y)}{\tau},      s \in[0,s_{11}] \label{eq:F1_c_x1} \\
		y_{11}(x,y,s)=&~\phi_1^{-1}(\phi_1(1)-s), \qquad \quad \,     s\in[0,s_{11}] \label{eq:F1_c_y1}\\
		x_{12}(x,y,s)=&~\frac{s}{\tau}, \qquad \qquad \qquad \qquad\,\,\,\, \,   s\in[0,s_{12}] \label{eq:F1_c_x2}  \\
		y_{12}(x,y,s)=&~\phi_1^{-1}(\tau x+\phi_1(y)-s), \,\,\,\,\,  s\in[0,s_{12}]\label{eq:F1_c_y2}
	\end{align}
\end{subequations}
where the argument $s \in [0,s_{1j} ]$, $j=1,\,2$, with $s_{1j}$ defined as 
\begin{subequations} \label{eq:s11_s12}
	\begin{align}
		s_{11}(x,y)=&~\phi_1(1)-\phi_1(y),\label{eq:s11}\\
		s_{12}(x,y)=&~\tau x.\label{eq:s12}
	\end{align}
\end{subequations}
(2) Characteristics of the $F_2(x,y)$:
For $(x,y)\in \mathcal{T}_2$, define the characteristic curves as follows, along which the system \eqref{ker_gen_2} becomes a family of ordinary differential equations,
\begin{subequations}
	\begin{align}
		x_{21}(x,y,s)=&~x+\frac{s}{\tau}-\frac{\phi_2(y)}{\tau},& s\in[0,s_{21}] \label{eq:F2_c_x1}\\
		y_{21}(x,y,s)=&~\phi_2^{-1}(s),& s\in[0,s_{21}]\label{eq:F2_c_y1}\\
		x_{22}(x,y,s)=&~ \frac{s}{\tau} , & s\in[0,s_{22}]\label{eq:F2_c_x2}\\
		y_{22}(x,y,s)=&~\phi_2^{-1}(\tau x+\phi_2(y)-s ),& s\in[0,s_{22}]  \label{eq:F2_c_y2}
	\end{align}
\end{subequations}
where the argument $s \in [0,s_{2j} ]$, $j=1,\,2$, with $s_{2j}$ defined as 
\begin{subequations}\label{eq:s21_s22}
	\begin{align}
		s_{21}(x,y)=&~\phi_2 (y), \label{eq:s21}\\
		s_{22}(x,y)=&~\tau x.\label{eq:s22}
	\end{align}
\end{subequations}
Integrating \eqref{ker_gen_1} along the characteristic curves defined by \eqref{eq:F1_c_x1}--\eqref{eq:F1_c_y1} and \eqref{eq:F1_c_x2}--\eqref{eq:F1_c_y2}, respectively, and considering boundary conditions \eqref{ker_gen_bud1} and \eqref{ker_gen_ini1}, respectively, the solution of $F_1(x,y)$ is yielded,\\
a) for  each $(x,y)\in \mathcal{T}^{F_1}_1$,
\begin{align}
	F_1(x,y)&= 
	G_{11}(x,y)+I_{11}[F](x,y),   \label{eq:F1_1}
\end{align}
b) for  each $(x,y)\in\mathcal{T}^{F_1}_2$,
\begin{align}
	F_1(x,y)&= 
	H_{1}(x,y)+G_{12}(x,y)+I_{12}[F](x,y),	 \label{eq:F1_2}
\end{align}
where $F$ is denoted 
$
	F=\begin{bmatrix}
		F_1  & F_2
	\end{bmatrix}^{\mathrm{T}}$,
and the functions $H_1$, $G_{1j}$ and $I_{1j}$, for $j=1,\,2$, will be defined later.

Similarly, integrating \eqref{ker_gen_2} along the characteristic curves defined by \eqref{eq:F2_c_x1}--\eqref{eq:F2_c_y1} and \eqref{eq:F2_c_x2}--\eqref{eq:F2_c_y2}, respectively, and considering boundary conditions \eqref{ker_gen_bud2} and \eqref{ker_gen_ini3}, respectively, the solution of $F_2  (x,y)$ is yielded, \\
a) for  each $(x,y)\in\mathcal{T}^{F_2}_1$,  
\begin{align}
	F_2=  Q_1(x) F_1(x ,0)+G_{21}(x,y) +I_{21}[F](x,y), \label{eq:F2_1} 
\end{align}
b) for each $(x,y)\in\mathcal{T}^{F_2}_2$, 
\begin{align}
	F_2&= H_{2}(x,y)+G_{22}(x,y)+I_{22}[F](x,y),\label{eq:F2_2}  
\end{align}
where the functions $H_i$, $G_{ij}$ and  $I_{ij}$, $i,\,j=1,\,2$  are defined as
\begin{subequations}
	\begin{align}
		&H_{i}(x,y) = h_i\big(y_{i2}(x,y,0)\big),\\
		&G_{ij}(x,y) = \! \int_{0}^{s_{ij}(x,y)}\!\! \! \!  \!\!     g_i\big(x_{ij}(x,y,s),y_{ij}(x,y,s)\big) \diff s,\\
		&I_{ij}[F](x,y)= \sum_{k=1}^{2} \int_{0}^{s_{ij}(x,y)}C_{ik}\big(x_{ij}(x,y,s),y_{ij}(x,y,s)\big) \cdot\nonumber\\
		& \quad\quad\quad\quad\quad  F_k\big(x_{ij}(x,y,s),y_{ij}(x,y,s)\big) \diff s.
	\end{align}
\end{subequations} 
Then, according to \eqref{eq:F1_1} and \eqref{eq:F1_2}, substituting $F_1(x,0)$ to \eqref{eq:F2_1}, the integral form of $F_2$ is rewritten as:\\
a) for  each $(x,y)\in\mathcal{T}^{F_1}_1\cap \mathcal{T}^{F_2}_1$,
\begin{align}
	F_2(x,y)=  &
	\psi_1(x,y)+Q_{11}  [F](x,y)+G_{21}(x,y)\nonumber\\
	&~+I_{21}[F](x,y), \label{eq:F2_new_1}
\end{align}
b) for  each $(x,y)\in\mathcal{T}^{F_1}_2\cap\mathcal{T}^{F_2}_1$,
\begin{align}
	F_2(x,y)=&~\psi_2(x,y)+Q_{12} [F](x,y)+G_{21}(x,y)\nonumber\\
	&~+I_{21}[F](x,y),\label{eq:F2_new_2}
\end{align}
c) for  each $(x,y)\in\mathcal{T}^{F_2}_2 $,
\begin{align}
	F_2(x,y)=H_{2}(x,y)+G_{22}(x,y)+I_{22}[F](x,y),\label{eq:F2_new_3}
\end{align}
where	
\begin{subequations}
	\begin{align}
		\psi_1 =&~q_1(x_{21}(x,y,0))G_{11}(x_{21}(x,y,0),0),\\
		\psi_2 =&~q_1(x_{22}(x,y,0))H_{1}(x_{22}(x,y,0),0)\nonumber\\
		&~+q_1(x_{22}(x,y,0))G_{12}(x_{22}(x,y,0),0),\\
		Q_{11}[F]  =&~q_1(x_{21}(x,y,0))I_{11}[F](x_{21}(x,y,0),0),\\
		Q_{12}[F]  =&~ q_1(x_{22}(x,y,0))I_{12}[F](x_{22}(x,y,0),0).
	\end{align} 
\end{subequations}

\subsection{ Convergence of the successive approximation series} \label{sec:convergence} 
Before proceeding, we first give the following lemmas.
\begin{lem} \label{lem:ineq_x_y}\normalfont
	If $(x,y) \in \mathcal{T}_2$ and $s \in [0,s_{ij} ]$, it holds that $(x_{ij}(x, y, s),y_{ij}(x, y, s)) \in\mathcal{T}_2 $ for $i,\,j =1,\,2$. Also, under the assumptions of Theorem \ref{theo:well}, $x_{ij}$, $y_{ij}$, and $s_{ij}$ are continuous in their domains of definition since they are defined as compositions of continuous functions. Moreover, the following inequalities are verified:
	\begin{subequations}
		\begin{align}
			x_{11}(x,y,s)&\leq x,\quad y_{11}(x,y,s)\geq y ,\\
			x_{12}(x,y,s)& \leq x,\quad y_{12}(x,y,s)  \geq y,\\
			x_{21}(x,y,s) &\leq x,\quad y_{21}(x,y,s) \leq y,\\
			x_{22}(x,y,s)& \leq x  ,\quad   y_{22}(x,y,s) \geq y.
		\end{align}
	\end{subequations}  
\end{lem}
The proof of this lemma is obvious, so it is omitted. This lemma is used in the proof of Lemma \ref{lem:x_int}.
\begin{lem}\label{lem:x_int}\normalfont For $n \in \mathbb{R}^+$, $i,\,j=1,2$, the following inequalities hold:
	\begin{align}
		\int_{0}^{s_{ij}(x,y)} x_{ij}^{n}(x,y,s)\diff s\leq   M_\varepsilon \frac{x^{n+1}}{n+1},
	\end{align}   
\end{lem}	 
 Similar to Lemma A.4 of \cite{coron2013lloca}, we can obtain this lemma by utilizing variable substitution and Lemma \ref{lem:ineq_x_y}. 

In this section, the successive approximation method is employed to prove that there is a solution of the system \eqref{eq:ker_gen} with the boundary condition \eqref{eq:ker_gen_bud}. We start by defining  the recursive equations for \eqref{eq:F1_1}--\eqref{eq:F1_2} and \eqref{eq:F2_new_1}--\eqref{eq:F2_new_3}  as follows, with $n=0,1,2,\cdots,+\infty$, 
\begin{subequations}\label{eq:F}
	\begin{align} 
		&F_1^{n+1}(x,y)=\begin{cases}\label{eq:F1}
			I_{11}[F^n](x,y), &\mathcal{T}^{F_1}_1,\\
			I_{12}[F^n](x,y), &\mathcal{T}^{F_1}_2,
		\end{cases}\\
		&F_2^{n+1}(x,y)=\\
		&  \begin{cases}
			Q_{11}(x) [F^n](x,y) +I_{21}[F^n](x,y), &\mathcal{T}^{F_1}_1\cap \mathcal{T}^{F_2}_1, \\
			Q_{12}(x)[F^n](x,y) +I_{21}[F^n](x,y), &\mathcal{T}^{F_1}_2\cap\mathcal{T}^{F_2}_1, \\
			I_{22}[F^n](x,y), &\mathcal{T}^{F_2}_2, \nonumber
		\end{cases}\label{eq:F2}
	\end{align}
\end{subequations}
with an initial function
\begin{subequations}
	\begin{align}
		F_1^0(x,y)&=\begin{cases}
			G_{11}(x,y) , &\mathcal{T}^{F_1}_1,\\
			H_{1}(x,y)+G_{12}(x,y) , &\mathcal{T}^{F_1}_2,
		\end{cases}\\
		F_2^0(x,y)&= \begin{cases}
			\psi_1(x,y) +G_{21}(x,y) , &\mathcal{T}^{F_1}_1\cap \mathcal{T}^{F_2}_1 ,\\
			\psi_2(x,y) +G_{21}(x,y) , &\mathcal{T}^{F_1}_2\cap\mathcal{T}^{F_2}_1 ,\\
			H_{2}(x,y)+G_{22}(x,y) , &\mathcal{T}^{F_2}_2. 
		\end{cases}
	\end{align}
\end{subequations}
It is easy to see that, if $ \sum_{n=0}^{+\infty} F^n(x,y)$ converges, then $F(x,y)$ can be alternatively written as $ F(x,y)=\sum_{n=0}^{+\infty} F^n(x,y)$,
where $F^n$ is denoted as
$
F^n=\begin{bmatrix}
	F_1^n & F_2^n 
\end{bmatrix}^{\mathrm{T}}$.

Next, we discuss the convergence of series by successive approximation. Define
\begin{subequations}
	\begin{align}   
		\bar C_{ij}=&~ \max_{(x,y)\in \mathcal{T}_{ij}}\left \{ \vert C_{ij}(x,y)\vert ,\, \vert\partial_x C_{ij}(x,y)\vert\right \},\\
		\bar q=&~ \max_{(x)\in [0,1]} \left \{ \vert q_1(x)\vert ,\, \vert\partial_xq_1(x)\vert \right\},\\
		M_0=&~\max\Bigg\{ \vert F_1^0 \vert, \, \vert F_2^0 \vert ,\, 2\sum_{k=1}^{2}\bar C_{ik}    \vert F_k^0\vert  ,\\
		&~\qquad \qquad \left(2\bar q \sum_{k=1}^{2}\bar C_{1k} +\sum_{k=1}^{2}\bar C_{2k} \right)   \vert F_k^0\vert \Bigg \},\\
		M_\varepsilon=&~\max \left\{  \tau , \, \frac{1}{\varepsilon_{1}(x)}, \, \frac{1}{\varepsilon_{2}(x)} \right\},\\
		\bar C=&~\max \left\{  3\bar{q} \sum_{k=1}^{2} \!\bar C_{1k} \!+\!2\sum_{k=1}^{2}\!\bar{C}_{2k},\, 3\sum_{k=1}^{2}\! \bar C_{ik}
		, \right\}\! ,
	\end{align}
\end{subequations}
with $i,j=1,2$. 
Assume that for  $n \in N^+$, \begin{align}\label{eq:ine_Fn} 
		\left\vert F^n(x,y) \right\vert \leq M_0 \frac{(\bar C  M_\varepsilon x)^n}{n!},
	\end{align}  
	holds, using induction method, we discuss whether $F^{n+1}(x,y)$ satisfies the above inequalities in different domain. From \eqref{eq:F} and \eqref{eq:ine_Fn}, using Lemma \ref{lem:x_int}, we have:\\
	(1) for $F_1^{n+1}$ on $\mathcal{T}_1^{F_1}$,
	\begin{align}
		&~\left \vert F_1^{n+1}(x,y)\right \vert =\left \vert I_{11}[F^n](x,y)\right \vert\nonumber\\
		\leq &~  \sum_{k=1}^{2}  \int_{0}^{s_{11}(x,y)}\Big \vert C_{1k}\big(x_{11}(x,y,s),y_{11}(x,y,s)\big)  \cdot \nonumber\\ &~\qquad\qquad   F_k^n\big(x_{11}(x,y,s),y_{11}(x,y,s)\big) \Big \vert \diff s   \nonumber\\
		\leq &~ M_0 \frac{(\bar C  M_\varepsilon )^n}{n!} \sum_{k=1}^{2} \bar C_{1k} \int_{0}^{s_{11}(x,y)}  x_{11}^n(x,y,s) \diff s \nonumber\\
		\leq&~ M_0 \frac{(\bar C  M_\varepsilon x)^{n+1}}{(n+1)!},
	\end{align}	
	(2) for $F_1^{n+1}$ on $\mathcal{T}_2^{F_1}$,
	\begin{align}
		&~ \left \vert F_1^{n+1}(x,y)\right \vert =\left \vert I_{12}[F^n](x,y)\right \vert\nonumber\\
		\leq&~  \sum_{k=1}^{2} \int_{0}^{s_{12}(x,y)}\Big \vert C_{1k}\big(x_{12}(x,y,s),y_{12}(x,y,s)\big)\cdot\nonumber\\
		&~\qquad \qquad F_k^n\big(x_{12}(x,y,s),y_{12}(x,y,s)\big) \Big \vert \diff s \nonumber\\
		\leq &~  M_0 \frac{(\bar C  M_\varepsilon )^n}{n!} \sum_{k=1}^{2} \bar C_{1k} \int_{0}^{s_{12}(x,y)}  x_{12}^n(x,y,s)  \diff s \nonumber\\
		\leq&~ M_0 \frac{(\bar C  M_\varepsilon x)^{n+1}}{(n+1)!},
	\end{align}	
	(3) for $F_2^{n+1}$ on $\mathcal{T}_1^{F_2} \cap \mathcal{T}_1^{F_1}$,
	\begin{align}
		&~ \left \vert F_2^{n+1}(x,y)\right \vert=\left \vert 
		Q_{11}(x) [F^n](x,y) +I_{21}[F^n](x,y)\right \vert\nonumber\\ 
		\leq&~ \left \vert q_1(x_{21}(x,y,0))I_{11}[F^n](x_{21}(x,y,0),0) \right \vert \nonumber\\
		&+\sum_{k=1}^{2} \int_{0}^{s_{21}(x,y)}\Big \vert C_{2k}\big(x_{21}(x,y,s),y_{21}(x,y,s)\big) \cdot\nonumber\\
		&~\qquad \qquad  F_k^n\big(x_{21}(x,y,s),y_{21}(x,y,s)\big) \Big \vert\diff s \nonumber\\
		\leq&~ \bar q \sum_{k=1}^{2} \int_{0}^{s_{21}(x_{21}(x,y,0),0)}\Big \vert C_{1k}\big(x_{21}(x,y,0),0\big)\cdot\nonumber\\
		&~\qquad \qquad\qquad \qquad F_k\big(x_{21}(x,y,0),0 \big) \Big \vert\diff s \nonumber\\&~+\sum_{k=1}^{2} \int_{0}^{s_{21}(x,y)}\Big \vert C_{2k}\big(x_{21}(x,y,s),y_{21}(x,y,s)\big)\cdot\nonumber\\
		&~\qquad \qquad \qquad F_k^n\big(x_{21}(x,y,s),y_{21}(x,y,s)\big)\Big \vert \diff s \nonumber\\
		\leq&~\bar{q} M_0 \frac{ \bar C ^n  M_\varepsilon ^{n+1}  } {n!}  \sum_{k=1}^{2} \bar C_{1k}x_{21}^{n+1}(x,y,0)\nonumber\\
		&~+ M_0 \frac{(\bar C^n  M_\varepsilon^{n+1}  )}{n!} \sum_{k=1}^{2}\bar{C}_{2k}x^{n+1}\nonumber\\
		\leq&~\left(\bar{q} \sum_{k=1}^{2} \bar C_{1k} +\sum_{k=1}^{2}\bar{C}_{2k}\right)M_0 \frac{ \bar C ^n  M_\varepsilon ^{n+1} x^{n+1} } {(n+1)!}  \nonumber\\
		\leq&~ M_0 \frac{(\bar C  M_\varepsilon x)^{n+1}}{(n+1)!},
	\end{align}	
	(4) for $F_2^{n+1}$ on $\mathcal{T}_1^{F_2} \cap \mathcal{T}_2^{F_1}$,
	\begin{align}
		&~\left\vert F_2^{n+1}(x,y)\right \vert =\left\vert Q_{12}(x)[F^n](x,y) +I_{21}[F^n](x,y)\right \vert \nonumber\\
		\leq&~\left\vert q_1(x_{22}(x,y,0))I_{12}[F](x_{22}(x,y,0),0)\right \vert\nonumber\\
		&~
		+\sum_{k=1}^{2} \int_{0}^{s_{21}(x,y)}\Big \vert C_{2k}\big(x_{21}(x,y,s),y_{21}(x,y,s)\big) \cdot\nonumber\\
		&~\qquad \qquad \qquad F_k^n\big(x_{21}(x,y,s),y_{21}(x,y,s)\big) \Big \vert\diff s \nonumber\\
		\leq &~ M_0 \frac{(\bar C  M_\varepsilon x)^{n+1}}{(n+1)!},
	\end{align}	
	(5) for $F_2^{n+1}$ on $\mathcal{T}_2^{F_2} $,
	\begin{align}
		&~\left \vert F_2^{n+1}(x,y)\right \vert=\left \vert I_{22}[F^n](x,y)\right \vert\nonumber\\
		&~\leq  \sum_{k=1}^{2} \int_{0}^{s_{22}(x,y)}\Big \vert C_{2k}\big(x_{22}(x,y,s),y_{22}(x,y,s)\big) \cdot \nonumber\\
		&~\qquad \qquad \qquad F_k^n\big(x_{22}(x,y,s),y_{22}(x,y,s)\big)\Big \vert \diff s \nonumber\\
		\leq &~  M_0 \frac{(\bar C  M_\varepsilon )^n}{n!} \sum_{k=1}^{2} \bar C_{2k} \int_{0}^{s_{22}(x,y)}  x_{22}^n(x,y,s)  \diff s \nonumber\\
		\leq&~ M_0 \frac{(\bar C  M_\varepsilon x)^{n+1}}{(n+1)!}.
	\end{align}
	Thus, by induction we have proved that \eqref{eq:ine_Fn} holds. Once the assumed value \eqref{eq:ine_Fn} is proved, it conclude that the series is bounded and converges uniformly,
	\begin{align}
		\left \vert \sum_{n=0}^{+\infty}   F^n(x,y)   \right\vert \leq  \sum_{n=0}^{+\infty}  M_0 \frac{(\bar C  M_\varepsilon x)^n}{n!} =M_0\rm{e}^{\bar C  M_\varepsilon x}.
	\end{align}
	Thus, it is proved that the system \eqref{eq:ker_gen} with the boundary condition \eqref{eq:ker_gen_bud} has a bounded solution $F(x,y)$.
	
	Next, we discuss the uniqueness of solution, assuming that two different solutions of \eqref{eq:ker_gen}--\eqref{eq:ker_gen_bud} are denoted by $F (x, y)$ and $F'(x, y)$. Define $\tilde F (x, y)=F (x, y)  - F'(x, y)$. Due to  the linearity of \eqref{eq:ker_gen}--\eqref{eq:ker_gen_bud}, $\tilde F(x,y)$ also verifies \eqref{eq:ker_gen}--\eqref{eq:ker_gen_bud}, with $q_1(x)=0$, $g_i(x,y) =0$ and $h_i(y) =0$, $i=1,\,2$. It is easy to get that $M_0 = 0$ for $\tilde F (x, y)$,  so we  get $\tilde F (x, y)=0$, which means $F (x, y) = F'(x, y)$.
	Thus, Theorem \ref{theo:well} is proved. 

\begin{figure*}[!ht]    
	\centering
	\begin{tabular}{@{}c@{}c@{}c@{}}
		\includegraphics[height=0.2\textwidth]{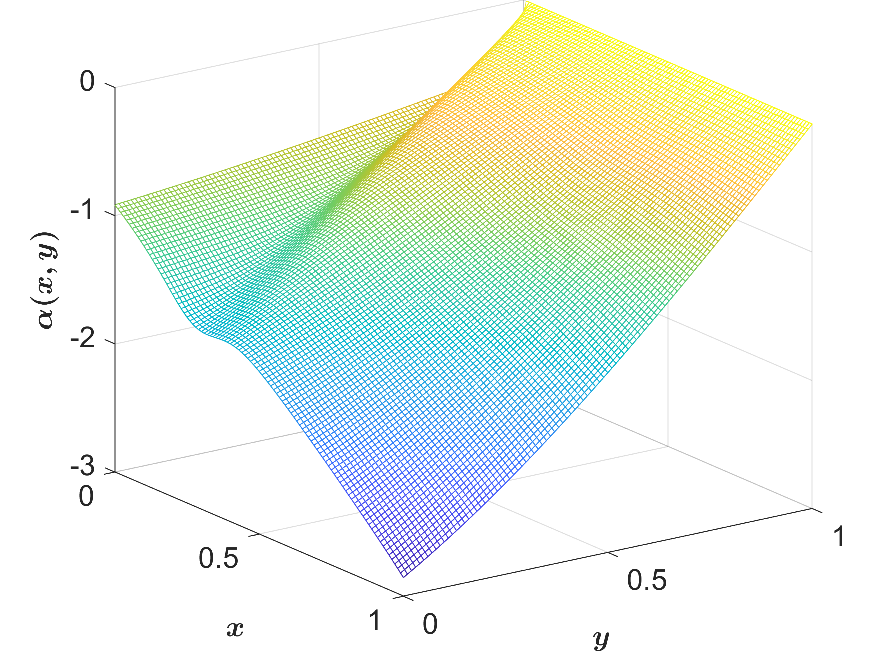}
		&\includegraphics[height=0.2\textwidth]{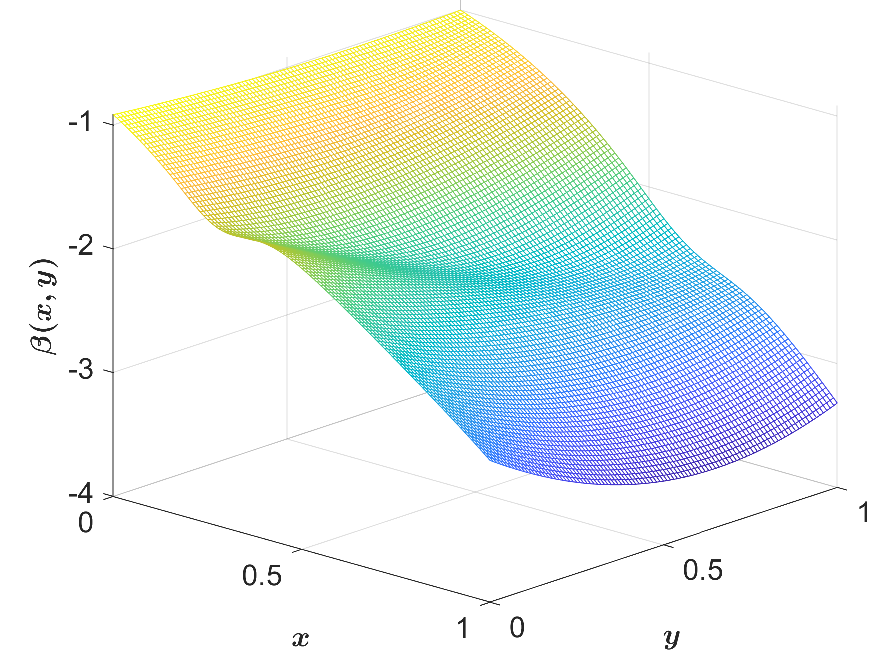}
		& \includegraphics[height=0.2\textwidth]{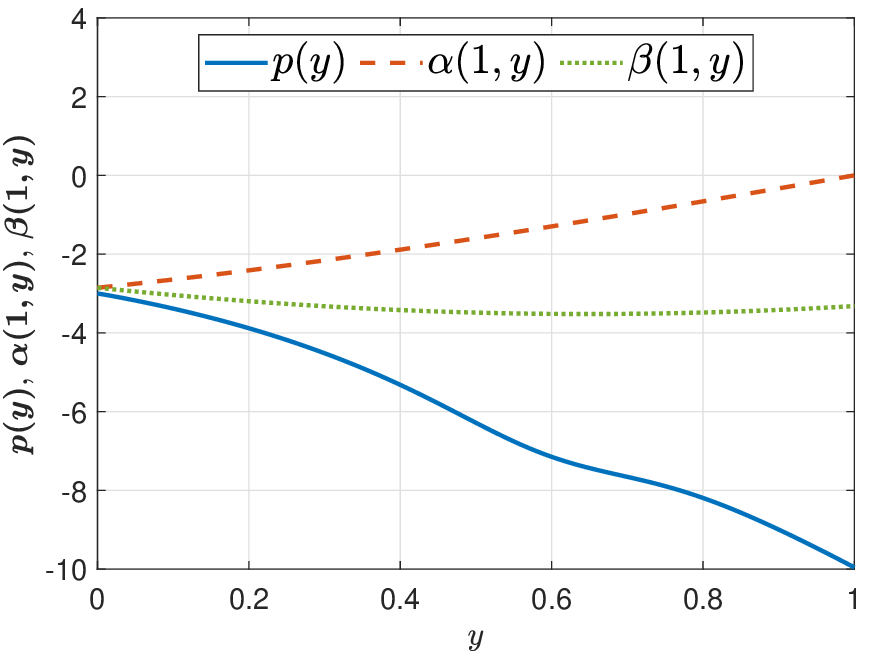} \\
		(a) $\alpha_1(x,y)$. &(b) $\alpha_2(x,y)$.  &(c)   Control gains $p(t)$, $\alpha_1(1,y)$ and $\alpha_2(1,y)$.
	\end{tabular}\\
	\caption{The numerical solution of kernel function.}
	\label{fig:kernel}
\end{figure*} 
\section{Finite-time Stability }\label{sec:stability}
In this section, we will establish the closed-loop stability. Before discussing the stability, we first give the invertibility of integration transformations;  then the equivalence between the target system and the cascaded system is constructed. Finally, the finite-time stability of the closed-loop system is obtained.

First, we present the inverse transformation $(w_1, w_2, z) \rightarrow (u_1, u_2, v)$, of transformation \eqref{eq:map} as follows:
\begin{subequations} \label{eq:map_inv}
	\begin{align}
		u_1(x ,t)=&~  w_1(x ,t) + \int_{0}^{x}  L^{11}(x ,y) w_1(y,t) \diff y\nonumber\\
		&~+\int_{0}^{x } L^{12}(x ,y) w_2(y,t) \diff y, \label{eq:inv_map1}\\
		u_2(x ,t)=&~   w_2(x ,t) + \int_{0}^{x }  L^{21}(x ,y) w_1(y,t) \diff y\nonumber\\
		&~+\int_{0}^{x } L^{22}(x ,y) w_2(y,t) \diff y, \label{eq:inv_map2}
		\\
		v(x ,t) 
		=&~z(x ,t)+\int_{0}^{x }\mu  (x -y) z(y,t) \diff y\label{eq:inv_map3}\\
		+ \int_{0}^{1}&  \beta_1  (x ,y) w_1(y,t) \diff y + \int_{0}^{1}  \beta_2(x ,y) w_2(y,t) \diff y.\nonumber 
	\end{align}
\end{subequations}
Well as the process in Section \ref{sec_sub:back}, using the integral by parts, we obtain the kernel functions $L^{ij}$, $i,j=1,2$ satisfy a set of hyperbolic PDEs given in \cite{coron2013lloca},
and the kernels $\beta_1(x,y)$, $\beta_2(x,y)$ satisfy:
\begin{subequations}\label{eq:beta}
	\begin{align}
		\!\frac{1}{\tau}\partial_x  \beta_1(x ,y) + \varepsilon_{1}(y) \partial_y \beta_1(x ,y) = & \! -\varepsilon_{1}'(y) \beta_1(x ,y) ,   \\  
		\!\frac{1}{\tau}\partial_x  \beta_2(x ,y) - \varepsilon_{2}(y) \partial_y \beta_2(x ,y)=&~\varepsilon_{2}'(y) \beta_2(x ,y), \\
		\beta_1(x ,1)=&~0,\\   
		\varepsilon_2(0)\beta_2(x ,0)=&~ q\varepsilon_{1}(0) \beta_1(x ,0),\\
		\beta_1(0,y)=&~L^{21}(1,y) ,\\
		\beta_2(0,y)=&~L^{22}(1,y) . 
	\end{align}
\end{subequations}
Then, according to the solution of $\beta_2(x,y)$, we get that
\begin{align}\label{eq:mu}
	\mu(x)=\tau \varepsilon_{2}(1) \beta_2(x
	,1).
\end{align} 
Second, the stability of target system \eqref{eq:target} is established by the following Lemma \ref{Prop:target}.
\begin{lem} \label{Prop:target} \normalfont
	Consider the system \eqref{eq:target} with initial conditions $w_{10}$, $w_{20}$ and $z_{0}  \in L^2([0,1])$,
	then the system is finite-time stable, that is,  its zero equilibrium is reached in finite-time $t = t_F$, where $
		t_F=\tau+\int_{0}^{1}\frac{1}{\varepsilon_{1}(y)}+\frac{1}{\varepsilon_{2}(y)} \diff y$.
\end{lem}
\begin{pf}
	\normalfont
	We give the explicit solution of \eqref{eq:target} by the method of characteristics as follows,
	\begin{subequations}
		\begin{align}
			w_1&=\left\{
			\begin{aligned}
				&	w_{10}\big(\phi_1^{-1}[\phi_1(x)-t]\big) ,&   t\leq  \phi_1(x),  \\
				&	qw_2(0,t-\phi_1(x))   ,  &t \textgreater \phi_1(x) ,
			\end{aligned}
			\right.\\
			w_2&=\left\{
			\begin{aligned}
				&	w_{20}\big(\phi_2^{-1}[t+\phi_2(x)]\big) , & t\leq \phi_2(1)-\phi_2(x) ,\\
				&	z (0,t+\phi_2(x)-\phi_2(1))   , & t\textgreater \phi_2(1)-\phi_2(x) ,
			\end{aligned}
			\right.\\
			z&=\left\{
			\begin{aligned}
				&	z_{0}\big(x+t/\tau\big) , & t\leq \tau(1-x) , \\
				&0   , & t\textgreater\tau,
			\end{aligned}
			\right.
		\end{align}
	\end{subequations}
	which concludes that the target system converges in finite time  $t_F$ defined in Lemma \ref{Prop:target}. Thus, this lemma has been proven.
\end{pf}
Through the above lemmas, the main result of this paper can be obtained as follows:
\begin{theorem}\normalfont \label{theo:mian}
	Consider the system \eqref{eq:coupled_system}  with the  control law \eqref{eq:U}. For any initial condition $u_{10}(x)$, $u_{20}(x)$, $v_{0}(x) \in L^2([0, 1])$, then the system is finite-time stable, that is, its zero equilibrium is reached in finite-time $t = t_F$,  where $t_F$ is given by Lemma \ref{Prop:target}.
\end{theorem}
Note that  the finite-time stability of system \eqref{eq:system} with control \eqref{eq:control}, which is equivalent to system \eqref{eq:coupled_system} and \eqref{eq:U}.
 
 \section{Robustness analysis} \label{sec:robust}	
  
 In this section, the delay-robustness is considered. Suppose that $\bar{\tau}=\tau+\Delta \tau>0$ is the uncertain delay which deviates from the true value $\tau$ by $\Delta \tau$. 
We map cascaded system \eqref{eq:coupled_system} by transformation \eqref{eq:map} to a new system which has same dynamics with \eqref{eq:main-coupled-w1}--\eqref{eq:target_v}, except the boundary condition: 
 \begin{align} \label{eq:bud-z}
 	z(1,t) =&~ U(t)-\int_{0}^{1}\mu  (1 -y) z(y,t) \diff y \\
 	-  \int_{0}^{1}  \beta_1 & (1 ,y) w_1(y,t) \diff y - \int_{0}^{1}  \beta_2(1,y) w_2(y,t) \diff y, \nonumber
 \end{align}
where  the inverse transformation \eqref{eq:map_inv} is used.
	Applying the method of characteristics on \eqref{eq:main-coupled-w1}--\eqref{eq:target_v} and \eqref{eq:bud-z}  yields,
	\begin{subequations}\label{eq:solution}
		\begin{align} 
			w_1(x,t) 
			=&~q z\left(1, t-\phi_1(x)-\phi_2(1)-\tau\right),\\ 
					w_2(x,t)  
			 =&~ z\left(1, t-\phi_2(1)+\phi_2(x)-\tau\right) ,\\
			 z(x,t)= &~z(1,t-\tau(1-x)).		
		\end{align}
	\end{subequations}
Substituting \eqref{eq:solution} into \eqref{eq:bud-z}, we have 
\begin{align}\label{eq:z1}  
		z(1,t)  =&~ U(t)-\int_{0}^{1}\mu  (y)z\left( 1, t-\tau y \right) \diff y  \nonumber\\
		- \int_{0}^{1}&   \beta_1  (1 ,y) q z\left(1, t-\phi_1(y)-\phi_2(1)-\tau\right) \diff y \nonumber\\
		- \int_{0}^{1} &   \beta_2(1,y) z\left(1, t-\phi_2(1)+\phi_2(y)-\tau\right) \diff y	.
\end{align}
	 \begin{remark} \normalfont
	 If one knows the exact delay $\tau$, then $z(1,t)=0$, which implies the \eqref{eq:main-coupled-w1}--\eqref{eq:target_v} and \eqref{eq:bud-z} is a finite-time stable system. However, the exact value of the delay is unknown, so the actual control $U(t)$ adopts   $\bar \tau$ which deviates from the true delay by a small value $\Delta \tau$. As a result, $z(1,t)$ captures the difference between the actual used control $U(t)$ and the desired control \eqref{eq:U}.
	 	\end{remark}


%
%

\begin{theorem}\label{theo:robust} \normalfont
  Consider the control law  adopting  $\bar\tau$: 
\begin{align}\label{eq:U_bar}   
U(t)= &~ \int_{0}^{1} \bar p (1-y) v(y,t) \diff y+\int_{0}^{1}  \bar \alpha_1 (1,y) u_1(y,t) \diff y\nonumber\\
&~  +\int_{0}^{1} \bar \alpha_2(1,y) u_2(y,t) \diff y, 
\end{align}
where  the functions $\bar p(\cdot)$ and $\bar \alpha_i(\cdot, \cdot)$, $i=1,2$, are  defined in the same way  as $ p(\cdot)$ and $ \alpha_i(\cdot, \cdot)$ substituting $\tau$  by $\bar \tau$ in \eqref{eq:alpha_main}--\eqref{eq:p}. There exist $\Delta \tau_0$ such that $\vert \bar \tau -   \tau\vert \leq \Delta \tau_0$,  the control \eqref{eq:U_bar}  robustly stabilizes the system \eqref{eq:coupled_system}.	 
\end{theorem}
\begin{pf}
	Suppose  $|\Delta \tau|=\vert\bar\tau-\tau\vert>0$.  The function $\bar \mu(\cdot)$ and $\bar \beta_i(\cdot, \cdot)$, $i=1,2$,  are defined by substituting $\bar \tau$ for $\tau$ in  \eqref{eq:beta}--\eqref{eq:mu}. Then, use the similar process as \eqref{eq:bud-z}--\eqref{eq:z1} for \eqref{eq:U_bar}, and substitute the result into \eqref{eq:z1}, which gives 
\begin{align}\label{eq:z111}
	 z(1,t)=& \int_{0}^{1}\bar\mu  (y) z\left( 1, t-\bar \tau y\right) \diff y  	-\int_{0}^{1} \mu  (y) z\left( 1, t-\tau y \right) \diff y  \nonumber\\
	+ \int_{0}^{1} &   \bar \beta_1  (1 ,y) q z\left(1, t-\phi_1(y)-\phi_2(1)-\bar\tau\right) \diff y \nonumber\\
	- \int_{0}^{1}&    \beta_1  (1 ,y) q z\left(1, t-\phi_1(y)-\phi_2(1)-\tau\right) \diff y \nonumber\\
	+ \int_{0}^{1}&  \bar \beta_2(1,y) z\left(1, t-\phi_2(1)+\phi_2(y)-\bar\tau\right) \diff y \nonumber\\ 
	- \int_{0}^{1}  &   \beta_2(1,y) z\left(1, t-\phi_2(1)+\phi_2(y)-\tau\right) \diff y.
\end{align}
Taking the Laplace transform of \eqref{eq:z111}, we get  characteristic equation  $P(s)=   1-\sum_{i=1}^{6}H_i(s)=0$, $i=1,2,\cdots,6$. The functions $H_i(s)$ are defined as,
\begin{align}
	H_1(s)=&\int_{0}^{1} \mu  (y)\mathrm{e}^{-   \tau y s}(\mathrm{e}^{-\Delta \tau y s}-1)  \diff y, \nonumber\\
	H_2(s)=&\int_{0}^{1} \Delta\mu(y)   \mathrm{e}^{- \bar \tau y s}\diff y ,\nonumber\\
	H_3(s)
	= & \int_{0}^{1} q \beta_1  (1 ,y) \mathrm{e}^{-(\phi_1(y)+\phi_2(1)+ \tau)s}(\mathrm{e}^{-\Delta \tau s}-1)       \diff y,\nonumber\\ 
 H_4(s)
	= &  \int_{0}^{1} q \Delta \beta_1  (1 ,y) \mathrm{e}^{-(\phi_1(y)+\phi_2(1)+ \bar \tau)s}      \diff y,\nonumber\\ 
	H_5(s)=&  \int_{0}^{1}   \beta_2  (1 ,y) \mathrm{e}^{-(\phi_2(1)-\phi_2(y)+ \tau)s}(\mathrm{e}^{-\Delta \tau s}-1)       \diff y,\nonumber\\
		 H_6(s)=&   \int_{0}^{1}  \Delta \beta_2  (1 ,y) \mathrm{e}^{-(\phi_2(1)-\phi_2(y)+\bar \tau)s}        \diff y,\nonumber
\end{align}
where $\Delta\mu(y) =\bar \mu(u)-\mu(y)$,    $\Delta\beta_1(1,y) =\bar \beta_1(1,y) -\beta_1(1,y)$ and $\Delta\beta_2(1,y) =\bar \beta_2(1,y)-\beta_2(1,y)$.
Due to Riemann-Lebesgues' lemma, we have
\begin{align}
	\forall\vert s\vert>M_0, \quad \int_{0}^{1} \mu  (y)\mathrm{e}^{-   \tau y s} \diff y \textless \frac{\xi}{6},
\end{align}
where $\xi\textless 1$. Then, consider the case of  $s\leq M_0$,   
\begin{align}\label{inequ:H1}
	|H_1(s)| \leq  \int_{0}^{1} |\mu (y)  |  \mathrm{e}^{-\tau ys} (\mathrm{e}^{|\Delta \tau| ys} -1)\ \diff y.
\end{align}
Since the continuity with respect to
$\tau$, we can  choose $\Delta \tau  \leq \Delta \tau_0$ small enough such that \eqref{inequ:H1} less than $\frac{\xi}{6}$. Hence, $\forall s\in \mathbb{C}$, $\forall \Delta \tau\le \Delta \tau_0$,   such that $\vert H_1(s)\vert \textless  \frac{\xi}{6}$.  Similarly, we have $\vert H_3(s)\vert \textless \frac{\xi}{6}$ and  $\vert H_5(s)\vert \textless \frac{\xi}{6}$. 
 
 Since the  continuity of $\mu$,   $\beta_1$ and $\beta_2$, we can choose $\Delta \tau  \leq \Delta \tau_0$ small enough such that  $\vert H_2(s)\vert \textless \frac{\xi}{6}$, $\vert H_4(s)\vert \textless \frac{\xi}{6}$ and $\vert H_6(s)\vert \textless \frac{\xi}{6}$. 
Consequently, we have
\begin{align}
 	\vert P(s)\vert  \geq 1 -\sum_{i=1}^{6}\vert H_i(s) \vert >0.  
\end{align}
It means that for  $0<\vert \Delta \tau \vert\leq \Delta \tau_0$,  $P(s)=0$  does  not have any roots on the complex right-open plane. As a result, \eqref{eq:z1} is asymptotically stable, which implies the system \eqref{eq:system} with control \eqref{eq:control} is robust stable due to the invertibility of the transformation \eqref{eq:map}.
\end{pf}
 
\begin{figure*}[!ht]    
	\centering
	\begin{tabular}{@{}c@{}c@{}c@{}}
		\includegraphics[height=0.2\textwidth]{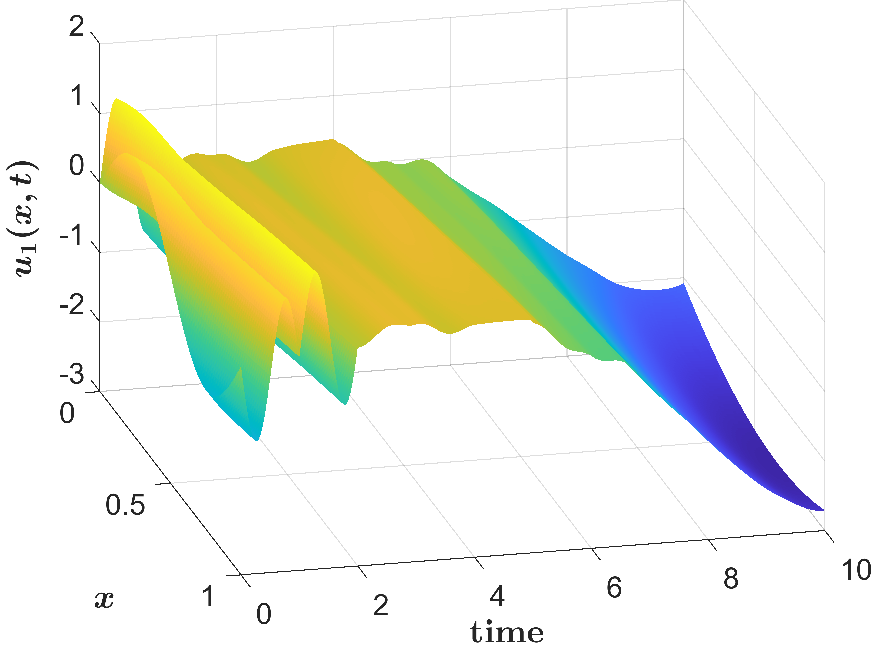}
		&\includegraphics[height=0.2\textwidth]{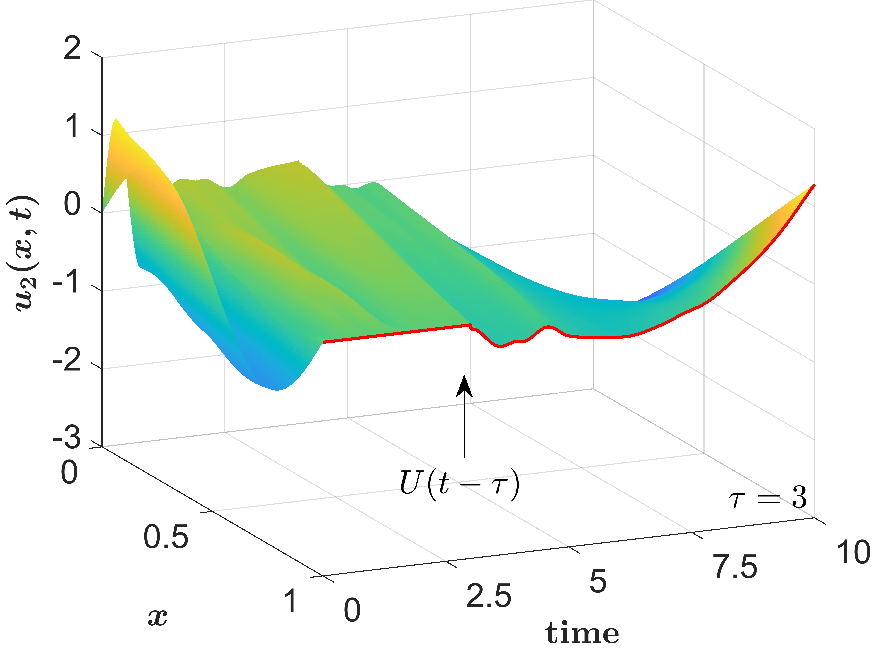}
		&\includegraphics[height=0.2\textwidth]{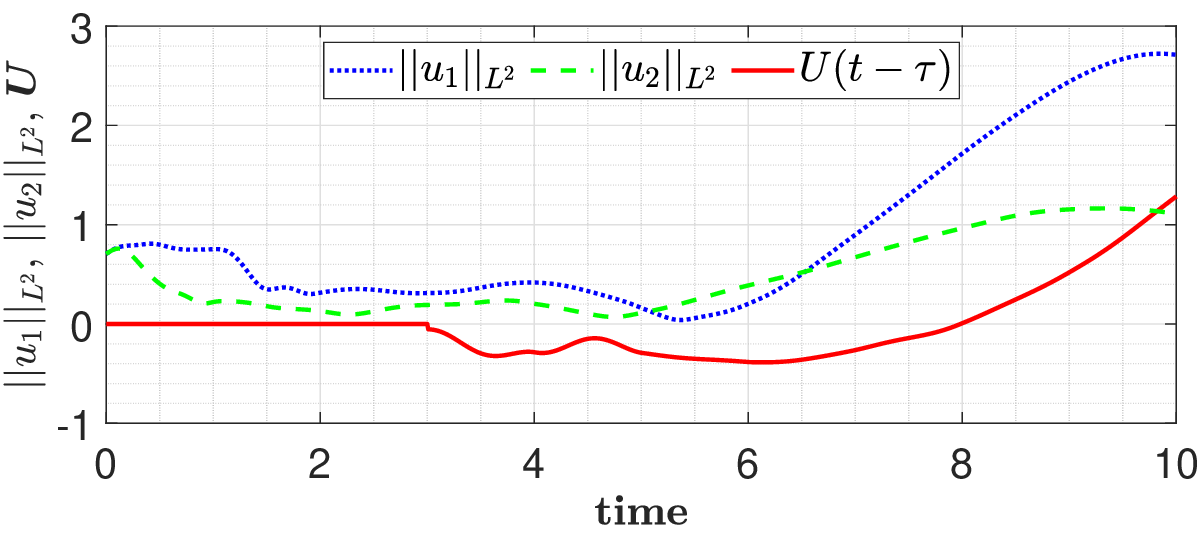} \\ 
		(a) The evolution of $u_1$. &(b) The evolution of $u_2$. &(c) The  $L_2$-norm and control effort. \\
		\includegraphics[height=0.2\textwidth]{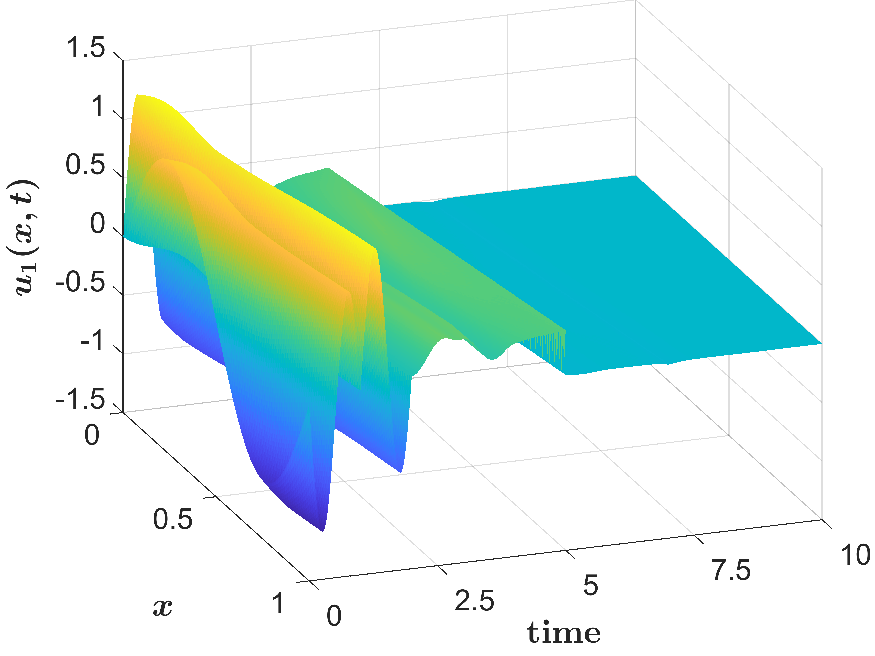} &\includegraphics[height=0.2\textwidth]{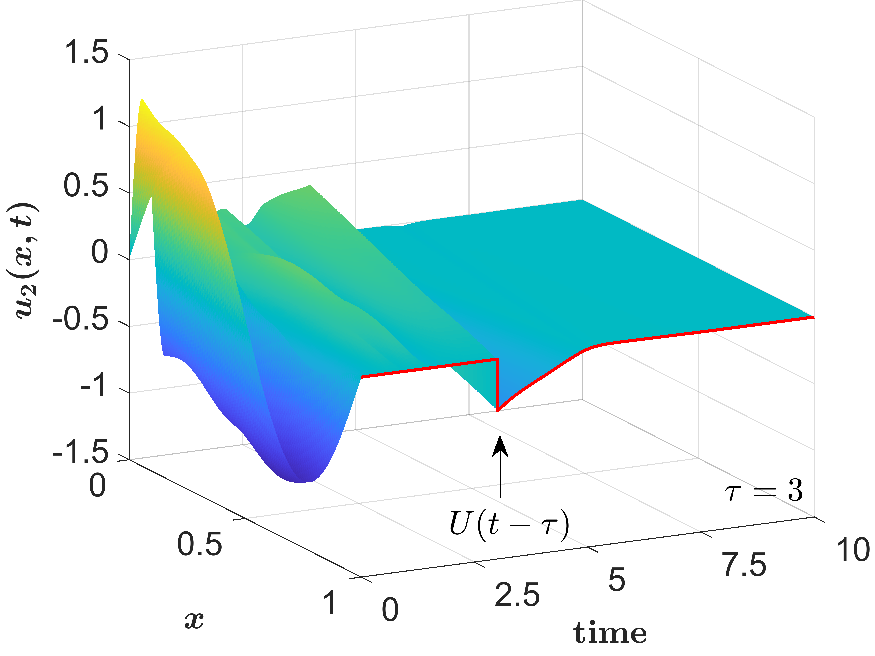}
		&\includegraphics[height=0.2\textwidth]{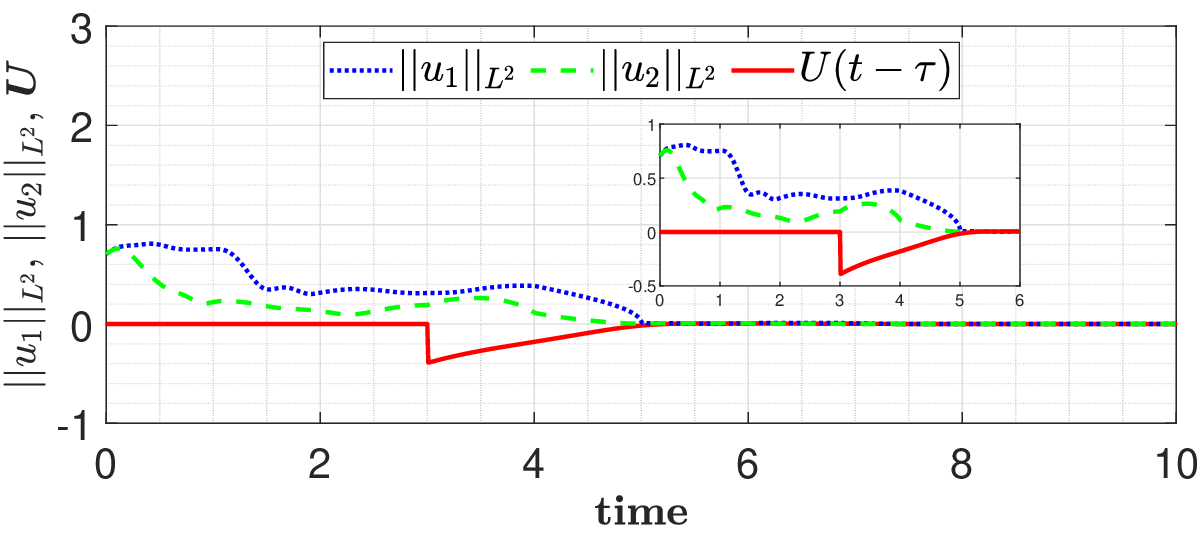}\\
		(d) The evolution of $u_1$. &(e) The evolution of $u_2$.&(f) The  $L_2$-norm and control effort. 
	\end{tabular}
	\caption{Top (a)--(b): The close-loop response of $2 \times 2$ hyperbolic system without delay compensation.   Bottom (c)--(d): The response under the delay compensator \eqref{eq:control}. }
	\label{fig:close-loop}
\end{figure*}
\begin{figure*}[!ht]    
	\centering
	\begin{tabular}{@{}c@{}c@{}c@{}}
		\includegraphics[height=0.2\textwidth]{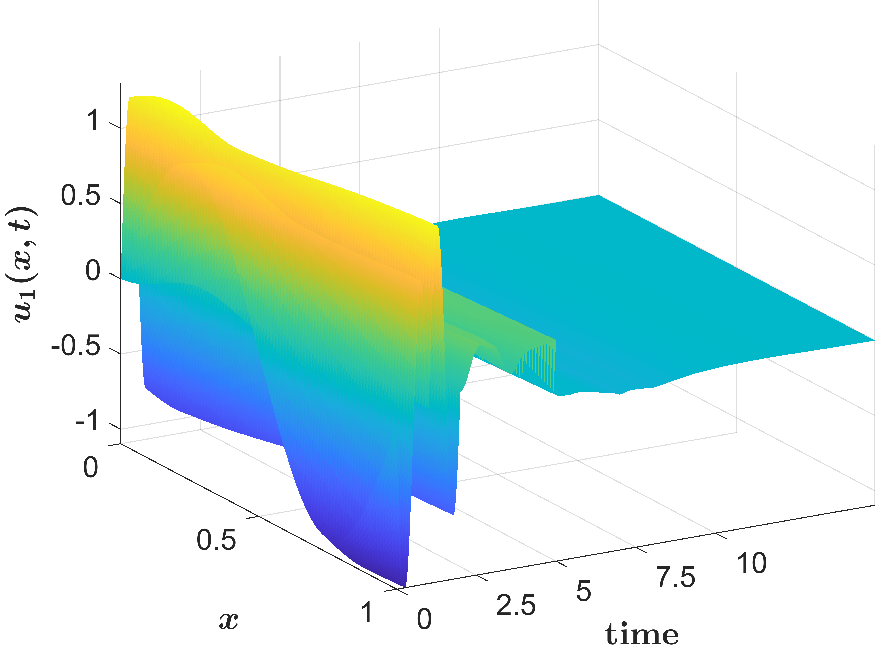} 
		&\includegraphics[height=0.2\textwidth]{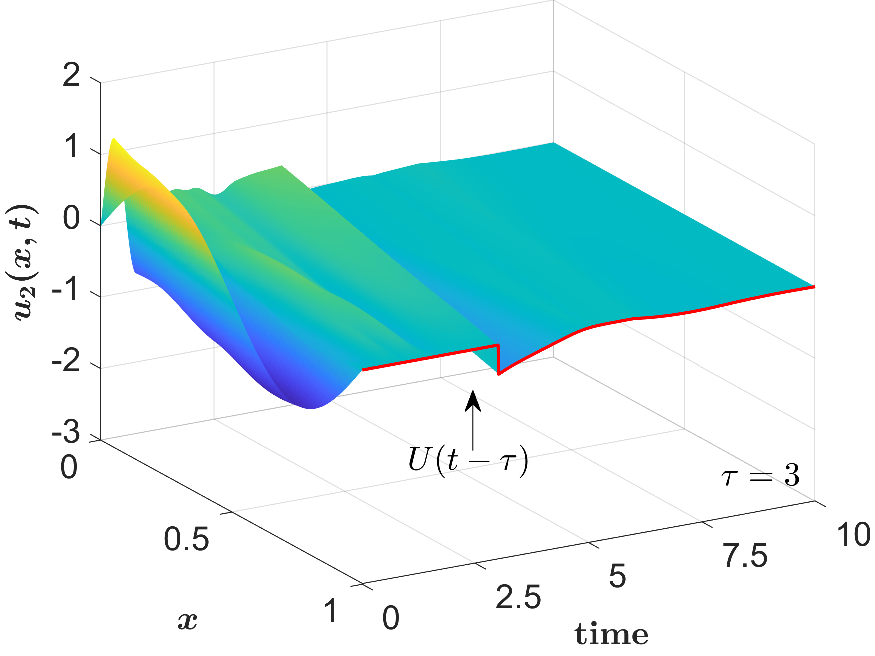}  &\includegraphics[height=0.2\textwidth]{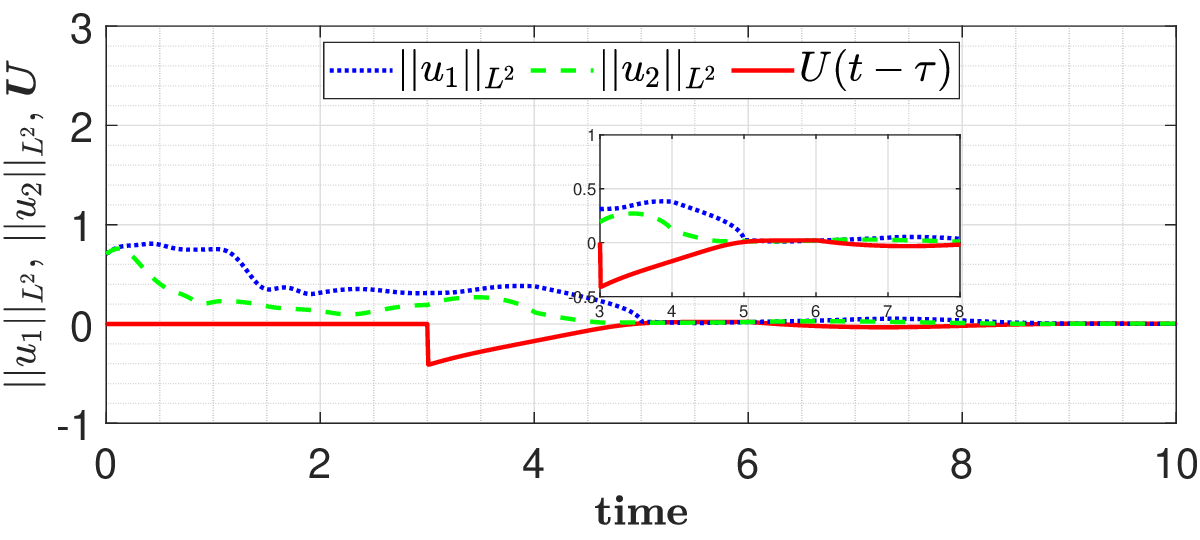}\\
		(a) The evolution of $u_1$. &(b) The evolution of $u_2$. &(c) The  $L_2$-norm and control effort.
	\end{tabular} 
	\caption{ The close-loop response under the delay compensator \eqref{eq:control}  in presence of delay uncertainty, where $\Delta \tau=0.2s$.} 
	\label{fig:robust}
\end{figure*}
\section{Simulation}\label{sec:simulation}
To illustrate Theorem \ref{theo:mian} and Theorem \ref{theo:robust}, in particular the role played the proposed compensator, the following examples are considered. The parameters in $2 \times 2$ first-order system \eqref{eq:system} are set as $\varepsilon_{1}=1$, $\varepsilon_{2}=1$, $c_1=1$, $c_2=1$, $q=1$ and $\tau=3s$. The initial state is set as
$u_{10}(x)=\sin(2\pi x)$, $u_{20}(x)=\sin(2\pi x)$.
The upwind scheme is employed in MATLAB, where time step
$h_t=0.01$ and space step $h_{x}=0.01$.  
The kernel function  $K^{ij}(x,y)$   can be solved explicitly using the technique in \cite{vazquez2014marcum}. Further, the kernel functions $\alpha_1(x,y)$, $\alpha_2(x,y)$ in \eqref{eq:alpha_main}--\eqref{eq:alpha_bud} and $p(y)$ in \eqref{eq:p} can be solved by using the finite difference method, and their numerical solutions are presented in Fig. \ref{fig:kernel}.
In Fig.\ref{fig:close-loop}(a)--(c), the delayed system evolution under the nominal controller from \cite{vazquez2011backstepping} indicates that a compensator is necessary.
It is easy to see that a closed-loop system without delay compensation cannot converge.
Meanwhile,  the system evolution under the proposed input delay compensator \eqref{eq:U} is shown in Fig. \ref{fig:close-loop}(d)--(e). Through simulation, it was proved that the system controlled by the compensation controller  converges to the zero equilibrium in $t_F=5s$, which agrees with Theorem \ref{theo:mian} perfectly.  Subsequently, Fig. \ref{fig:robust} depicts the evolution of the closed-loop system in presence of an uncertain input delay. The input delay is $\tau=3s$ in the actual system, but the parameter used in the compensator with delay variation  $\Delta \tau=0.2s$. As expected in Theorem \ref{theo:robust}, the closed-loop system is robustly stable.


\section{Conclusion}\label{sec:conclusion}
In this paper, we propose a backstepping-based input-delay-compensated boundary controller to stabilize a $2 \times 2$   hyperbolic system subject to an uncertain  input delay. The well-posedness of the kernel function is explored using the method of characteristics and successive approximation method. Therefore, under the action of the proposed compensator, the finite-time stability of the hyperbolic system are guaranteed.  Moreover, the robustness of the proposed controller is analyzed using the equivalent neutral equation. Future research issues are to consider compensation control design with non-constant delays, such as spatial- or time-varying delays, or to consider the time-delay compensation problem for more general $m+n$ hyperbolic system.




\bibliographystyle{elsarticle-num} 
\bibliography{First_order1.bib}           
\end{document}